\documentclass[accepted]{uai2023} 

\usepackage[british]{babel}

\usepackage[utf8]{inputenc}
\usepackage{philip}

\newcommand\xqed[1]{%
  \leavevmode\unskip\penalty9999 \hbox{}\nobreak\hfill
  \quad\hbox{#1}}
\newcommand\exqed{\xqed{$\triangle$}}

\title{Correcting for Selection Bias and Missing Response \\ in Regression using Privileged Information}

\author[1,2]{\href{mailto:philip.boeken@gmail.com?subject=Your UAI 2023 paper}{Philip~Boeken}}
\author[1]{Noud~de~Kroon}
\author[2]{Mathijs~de~Jong}
\author[1]{Joris~M.~Mooij}
\author[2]{Onno~Zoeter}
\affil[1]{%
    Korteweg-de Vries Institute for Mathematics\\
    University of Amsterdam\\
    The Netherlands
}
\affil[2]{%
    Booking.com\\
    The Netherlands
}

\begin{document}

\maketitle

\begin{abstract}
    When estimating a regression model, we might have data where some labels are missing, or our data might be biased by a selection mechanism.
    When the response or selection mechanism is \emph{ignorable} (i.e., independent of the response variable given the features) one can use off-the-shelf regression methods; in the \emph{nonignorable} case one typically has to adjust for bias. We observe that \emph{privileged information} (i.e.\ information that is only available during training) might render a nonignorable selection mechanism ignorable, and we refer to this scenario as \emph{Privilegedly Missing at Random} (PMAR). We propose a novel imputation-based regression method, named \emph{repeated regression}, that is suitable for PMAR. We also consider an importance weighted regression method, and a doubly robust combination of the two. The proposed methods are easy to implement with most popular out-of-the-box regression algorithms. We empirically assess the performance of the proposed methods with extensive simulated experiments and on a synthetically augmented real-world dataset. We conclude that repeated regression can appropriately correct for bias, and can have considerable advantage over weighted regression, especially when extrapolating to regions of the feature space where response is never observed.
\end{abstract}


\section{Introduction}
Regression is a primary technique in data science, machine learning and statistics. When presented with data $(X_1, Y_1), ..., (X_n, Y_n)$ sampled from the distribution $\PP(X, Y)$, the goal is to find the conditional expectation of $Y$ given $X$, i.e.\ estimate the function $\mu(x) = \EE[Y|X=x]$. Practitioners are often presented with either incomplete data (where some values are missing) or data that is not representative of the population (drawn from some other $\tilde{\PP}(X, Y)$ with $\tilde{\PP} \neq \PP$), and hence have to correct for the bias that is present in their training data. The discrepancy between $\tilde{\PP}$ and $\PP$ can for example be induced by a selection mechanism.

Consider the situation where we want to predict whether a loan applicant will default. Let $X$ be the digital record of a person applying for a loan, $Y$ whether the borrower defaults on the loan, $\hat{Y}$ a current algorithm’s prediction of default, and $Z$ (costly) expert advice on whether the applicant will default, which correlates with $Y$ through information that is unavailable to us. We want to reject any applicants that will default on the loan, so we have the issue of the loan $S$, as a weighted combination of $\hat{Y}$ and $Z$. Only when we issue the loan, $Y$ will be observed. If we want to re-train our current prediction model $\hat{Y} = \hat{\EE}[Y|X]$ we have to take into account the bias induced by $Z$ during training of the model without explicitly adding $Z$ to the covariates $X$. In this work, we demonstrate how one can incorporate such \emph{privileged information} in a regression model, to correct for any bias that it induces in the data generating process.

An important line of work on selection bias is by \cite{pearl2012solution} and \cite{bareinboim2014recovering}, who consider the problem of estimating $\PP(Y|X)$ from a potentially biased dataset by leveraging knowledge of the underlying causal graph. They derive an expression of $\PP(Y|X)$ as an integral of quantities that can be estimated from readily available biased and additional (`external') unbiased data. Although identification of a conditional distribution might be useful in certain scenarios, this does not tell the practitioner how to estimate a regression model, especially when dealing with continuous variables. In this work we address this problem, keeping in mind the applicability of the proposed methods.

Missingness problems are often characterised by the presence of certain conditional independencies in the data. Typically these independencies are untestable, making it unclear whether these conditional independence assumptions are appropriate for the data at hand. Drawing a causal graph of the data generating process is often helpful for gaining a better understanding of the problem. Moreover, the causal graph identifies conditional independencies in the data, and thereby allows the practitioner to motivate these conditional independence assumptions. As posed by \cite{rubin1976inference}:

\textit{``The inescapable conclusion seems to be that when dealing with real data, the practising statistician should explicitly consider the process that causes missing data far more often than he does. However, to do so, he needs models for this process and these have not received much attention in the statistical literature.''}

Contrary to what is often seen in the literature we do not need a full probabilistic model of the missingness mechanism, as the identification of certain conditional independencies can already be sufficient for estimating a regression model.

Our contributions are as follows. First, this paper serves a pedagogical purpose by reviewing literature of missingness and selection bias, as addressed in Section \ref{sec:missingness}. We address what the consequences of certain conditional independencies are for the practitioner. While doing so, we motivate the use of \emph{privileged information} \citep{vapnik2009new} when missingness or selection is nonignorable, and introduce the \emph{Privilegedly Missing at Random} setting (PMAR).
Secondly, in Section \ref{sec:methods} we formulate a novel imputation-based estimator for the PMAR setting, which we refer to as \emph{repeated regression}. We point out that importance weighted regression is also suitable for the PMAR setting, and combine the two methods into a doubly robust estimator.
These estimators are formulated such that they are easy to implement using out-of-the box regression algorithms.
We then address the intricacies of evaluation under selection bias and missingness in Section \ref{sec:evaluation}. We warn the practitioner that, to our knowledge, there is no appropriate way of evaluating a regression method on a finite, biased dataset without relying on auxiliary models.
Lastly, in Section \ref{sec:experiments} we assess the performance of the formulated methods on extensive simulated experiments and on synthetically augmented real-world data. We observe that repeated regression has considerable advantage over importance weighting, especially when extrapolating.

\section{Missing data, selection bias and privileged information}\label{sec:missingness}
\subsection{Missing data mechanisms}
A very general framework for handling missing data is proposed by \cite{rubin1976inference}. When modelling the distribution of a set of random variables $X_1, ..., X_m$, this framework takes into account potential missing values of any of the $X_i$ by considering a random vector of \emph{response indicators}: binary variables $S_1, ..., S_m$, where $S_i = 1$ indicates that variable $X_i$ is observed. The process that determines whether we observe $X_i$ can be explicitly modelled as $\PP(S_i | X_1, ..., X_m)$. We then define the \emph{observed} random vector $X_1^o, ..., X_m^o$ where $X_i^o = X_i$ if $S_i=1$ and $X_i = ?$ if $S_i = 0$, where `?' denotes the value that is missing in the dataset. We focus on a specific missingness problem where we have covariates $X$, response variable $Y$, and where $S$ is a response indicator for $Y$, so only $Y$ can have missing values.

\cite{rubin1976inference} proposed the following classification of missing data mechanisms, which approximately categorises the difficulty of many inference problems.
\begin{definition}{\citep{rubin1976inference}}
    Given variables $X, Y$, and response indicator $S$ for $Y$, we say that $Y$ is
    \begin{itemize}
        \item Missing Completely at Random (MCAR) if the missingness mechanism is independent of all other observed variables, i.e.\ $S\Indep X, Y$;
        \item Missing at Random (MAR) if the missingness mechanism is independent of the missing variable given all other fully observed variables, i.e.\ $Y \Indep S \given X$;
        \item Missing Not at Random (MNAR) if it is neither MCAR nor MAR.
    \end{itemize}
\end{definition}

A large body of literature has been written about inference under missingness. Important works are the EM algorithm \citep{dempster1977maximum}, Nobel prize winning work by \cite{heckman1979sample} on correcting for selection bias in linear regression, \cite{rosenbaum1984reducing} which laid the foundations of propensity score based methods in causal inference, and a series of generalised estimating equations (GEE) based methods for inference under missingness \citep{robins1994estimation, robins1995semiparametric,rotnitzky1998semiparametric,scharfstein1999adjusting}.
Graphical modelling of missingness mechanisms and related independence testing and identification problems have been investigated in \cite{daniel2012using,thoemmes2015graphical,nabi2020full,mohan2021graphical,goel2021importance}. For an overview of the field, see \cite{little2019statistical}. In this work we focus on correcting for bias in regression using privileged data, which is not treated in the aforementioned literature.

\paragraph{S-recoverability}
As proposed by \cite{pearl2012solution} and further developed in \cite{bareinboim2014recovering}, \emph{s-recoverability} is a method to deal with selection bias.
Consider again the case of modelling the distribution of $X_1, ..., X_m$. Instead of having one response indicator for each $X_i$ as in the missingness framework, \cite{pearl2012solution} considers one selection variable $S$ where $S=1$ indicates that all variables are observed (\emph{selected}), and $S=0$ indicates that no variable is observed, i.e. there is no row for this observation in our dataset. Any data that we observe is drawn from the distribution $\PP(X_1, ..., X_m | S=1)$.

\cite{bareinboim2014recovering} prove that for discrete variables and under certain positivity assumptions, $\PP(Y|X)$ is recoverable from $\PP(X, Y|S=1)$ if and only if $Y\Indep S \given X$. The `if' part is straightforward, since the conditional independence implies $\PP(Y|X) = \PP(Y|X, S=1)$, and the right-hand-side can be estimated from the data.
When this conditional independence is not satisfied, \cite{bareinboim2014recovering} consider joint measurement of $X, Y$ and some other variable $Z$ in the biased dataset (so availability of $\PP(X, Y, Z | S=1)$), satisfying the conditional independence $Y\Indep S\given X, Z$. Additionally, they assume availability of unbiased measurements of $(X, Z)$, i.e.\ sampled from $\PP(X, Z)$, which they refer to as `external data'. For discrete $Z$, the quantity $\PP(Y|X)$ is then identified with
\begin{align}\label{eqn:s-recoverability}
    \begin{split}
        \PP(Y|X) &= \sum_{z} \PP(Y|X, Z=z) \PP(Z=z|X) \\
        &= \sum_{z} \PP(Y|X, Z=z, S=1) \PP(Z=z|X)
    \end{split}
\end{align}
provided we have $\PP(S=1 | X, Z) > 0, \PP(X, Z)$-almost surely, which is for example satisfied when $\supp[\PP(X, Z | S=1)] = \supp[\PP(X, Z)]$.\footnote{We let $\supp[\PP]$ denote the support of $\PP$.}
One can straightforwardly replace the sum with an integral when $Z$ is continuous. However, when the domain of $Z$ is countably infinite or continuous, estimating this quantity is not straightforward. The repeated regression estimator proposed in Section \ref{subsec:rr} is a solution to this problem.

Our proposed methods can be applied to the selection bias and missingness settings. In either case, we require availability of i.i.d.\ observations of $(X, Y, Z) \sim \PP(X, Y, Z | S=1)$ whose index set we denote with $\Scal$, and i.i.d.\ observations of $(X, Z) \sim \PP(X, Z)$ whose index set we denote with $\Dcal$. In the missingness setting, both samples are readily available and we have $\Scal \subseteq \Dcal$. In the selection bias setting the sample $\Dcal$ consists of `external' data and we typically have $\Scal \cap \Dcal = \emptyset$. A schematic display of the assumed available data under missingness and selection bias is provided in the supplements. For the methods of Sections \ref{subsec:ipw} and \ref{subsec:dr} we additionally require knowledge of the selection probability $\PP(S=1 | X, Z)$, which is directly estimable in the missingness setting, but which has to be assumed under selection bias. Throughout this paper, any distinctions between missingness and selection bias will be pointed out when necessary. Otherwise, either setting can be assumed.

\subsection{Regression under different selection mechanisms}
Suppose we are interested in estimating $\EE[Y|X]$ with continuous or discrete (ordinal) $Y$, and arbitrary $X$.\footnote{Note that this setting includes binary classification.} We might be confronted with a dataset with missing values of $Y$, or we might suspect that some selection mechanism is in play which makes certain $(X, Y)$ pairs unobserved. In this section, we investigate whether there is need for any bias correction. As this investigation is based on conditional independence assumptions, we first elaborate how such assumptions can be motivated.

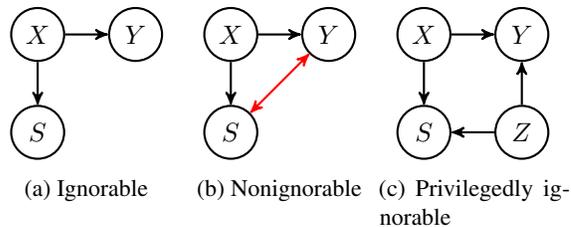
\begin{figure}
    \centering
    \begin{subfigure}[t]{.3\linewidth}
        \centering
        \begin{tikzpicture}
            \node[var] (Y) at (1.3, 0) {$Y$};
            \node[var] (X) at (0, 0) {$X$};
            \node[var] (S) at (0, -1.3) {$S$};
            \draw[arr] (X) to (Y);
            \draw[arr] (X) to (S);
        \end{tikzpicture}
        \caption{Ignorable}
        \label{fig:ignorable}
    \end{subfigure}
    \begin{subfigure}[t]{.3\linewidth}
        \centering
        \begin{tikzpicture}
            \node[var] (Y) at (1.3, 0) {$Y$};
            \node[var] (X) at (0, 0) {$X$};
            \node[var] (S) at (0, -1.3) {$S$};
            \draw[arr] (X) to (Y);
            \draw[arr] (X) to (S);
            \draw[biarr, color=red] (Y) to (S);
        \end{tikzpicture}
        \caption{Nonignorable}
        \label{fig:non_ignorable}
    \end{subfigure}
    \begin{subfigure}[t]{.3\linewidth}
        \centering
        \begin{tikzpicture}
            \node[var] (Y) at (1.3, 0) {$Y$};
            \node[var] (Z) at (1.3, -1.3) {$Z$};
            \node[var] (X) at (0, 0) {$X$};
            \node[var] (S) at (0, -1.3) {$S$};
            \draw[arr] (X) to (Y);
            \draw[arr] (X) to (S);
            \draw[arr] (Z) to (Y);
            \draw[arr] (Z) to (S);
        \end{tikzpicture}
        \caption{Privilegedly ignorable}
        \label{fig:cond_ignorable}
    \end{subfigure}
    \caption{Examples of missingness or selection bias settings, where $S$ indicates whether $Y$ is observed or not.}
    \label{fig:missingness}
\end{figure}

\paragraph{Causal modelling} Missingness mechanisms can be characterised by independencies. This is everything we need: all proposed methods will only require certain conditional independencies in the data, and no causal assumptions. However, the conditional independencies that are assumed are typically untestable. For example, the independence $Y\Indep S$ is not testable, as we have not observed $Y$ for $S=0$. To motivate such an independence assumption, one could model the data generating process with a graphical causal model, and infer from d-separations in the graph that there must be certain independencies in the data \citep{pearl2009causality}.
In the missingness setting where only $Y$ can be missing (with indicator $S$), we can draw a simplified graph that discards the $Y^o$ variable, as in Figure \ref{fig:ignorable}. In the graphical framework of s-recoverability, $S$ is implicitly required to be a sink node (i.e.\ a node without children) \citep{bareinboim2014recovering}. We do not adopt this convention.

\paragraph{Ignorable missingness/selection bias} When the data is MCAR or MAR, we have the conditional independence $Y\Indep S \given X$ and thus
\begin{equation}\label{eqn:id_ign}
    \EE[Y|X] = \EE[Y|X, S=1],
\end{equation}
in which case we call the missingness/selection mechanism \emph{ignorable} for estimating $\EE[Y|X]$. An example of an ignorable (MAR) missingness or selection mechanism is given in Figure \ref{fig:ignorable}. When this mechanism is ignorable, any correctly specified model could directly be learned as \emph{Empirical Risk Minimization} (ERM) is consistent \citep{sugiyama2007covariate}, and hence, there is no need for bias correction. Note that the right-hand side of equation (\ref{eqn:id_ign}) is only defined for $X \in \supp[\PP(X | S=1)]$, and one has to be mindful when extrapolating to values $X\in\supp[\PP(X)] \setminus \supp[\PP(X | S=1)]$ \citep{martius2016extrapolation}.

Although in this case the missingness or selection mechanism can in principle be ignored, one should beware that under model misspecification, ERM is not consistent anymore. Also, efficiency of the estimation procedure can be affected by MCAR and MAR.
\cite{zadrozny2004learning} and \cite{weifan2005improved} investigate the performance of multiple popular classification algorithms under covariate shift.
\cite{robins1995semiparametric} investigate the asymptotic efficiency of semi-parametric estimation of the regression function via Generalised Estimating Equations (GEE). Note that in the missingness setting, we can use the values of $X$ where $Y$ is unobserved as additional (unlabelled) input of a \emph{semi-supervised learning} algorithm. When doing \emph{anticausal} prediction (i.e., when $Y$ is the cause of $X$), these additional samples could benefit performance \citep{scholkopf2012causal}.

\paragraph{Privilegedly ignorable missingness/selection bias}
It might be the case that the missingness or selection mechanism is nonignorable. For example, we may suspect that $S$ and $Y$ are confounded by some variable that is not contained in $X$, as depicted in Figure \ref{fig:non_ignorable}. In this case we typically have $\EE[Y|X, S=1] \neq \EE[Y|X]$, in which case naive regression on the available data yields a biased model. However, we might be able to identify the latent confounder, or more generally, any variable (e.g.\ mediator or confounder of $Y$ and $S$) that induces a dependence between $Y$ and $S$ after conditioning on $X$.

More specifically, we are looking for a set of variables $Z$ such that $Y\Indep S \given X, Z$, such as depicted in Figure \ref{fig:cond_ignorable}. Having identified $Z$, one might include it into the features of the model and deploy the unbiased model $\hat{\EE}[Y|X, Z, S=1]$. However, in practice it might be costly or impossible to measure $Z$ during deployment, e.g.\ when $Z$ can only be measured after making a prediction with the regression model. In these settings, it is undesirable to incorporate $Z$ into the features of the deployed model. When $Z$ is only available during training it is referred to as \emph{privileged information}, following \cite{vapnik2009new}. The specific setting that we consider in this work is formally defined as follows:
\begin{definition}[PMAR]\label{def:pmar}
    Given features $X$ and label $Y$, response indicator $S$ for $Y$, and privilegedly observed variable $Z$, we say that $Y$ is
    \emph{Privilegedly Missing at Random} (PMAR)
    if the response indicator is independent of the missing label given all other fully or privilegedly observed variables, i.e.\ $Y \Indep S \given X, Z$.
\end{definition}
When the data used for regression is PMAR, we refer to the missingness/selection mechanism as \emph{privilegedly ignorable}.

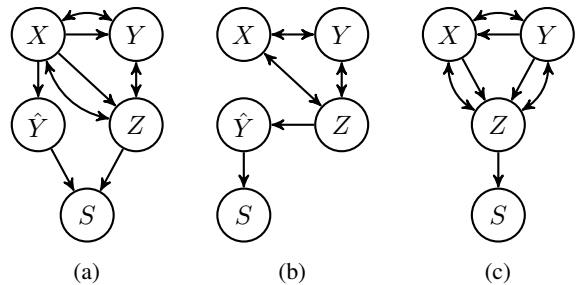
\begin{figure}[tb]
    \centering
    \begin{subfigure}[t]{0.32\linewidth}
        \centering
        \begin{tikzpicture}
            \node[var] (X) at (0, 0) {$X$};
            \node[var] (Yhat) at (0, -1.2) {$\hat{Y}$};
            \node[var] (Y) at (1.3, 0) {$Y$};
            \node[var] (Z) at (1.3, -1.2) {$Z$};
            \node[var] (S) at (0.65, -2.4) {$S$};
            \draw[biarr, bend left] (X) to (Y);
            \draw[arr] (X) to (Y);
            \draw[biarr, bend right] (X) to (Z);
            \draw[arr] (X) to (Z);
            \draw[arr] (X) to (Yhat);
            \draw[arr] (Yhat) to (S);
            \draw[arr] (Z) to (S);
            \draw[biarr] (Z) to (Y);
        \end{tikzpicture}
        \caption{}
        \label{fig:pmar:selection:1}
    \end{subfigure}
    \begin{subfigure}[t]{0.32\linewidth}
        \centering
        \begin{tikzpicture}
            \node[var] (X) at (0, 0) {$X$};
            \node[var] (Z) at (1.3, -1.2) {$Z$};
            \node[var] (Yhat) at (0, -1.2) {$\hat{Y}$};
            \node[var] (Y) at (1.3, 0) {$Y$};
            \node[var] (S) at (0, -2.4) {$S$};
            \draw[biarr] (Z) to (Y);
            \draw[biarr] (X) to (Y);
            \draw[biarr] (X) to (Z);
            \draw[arr] (Z) to (Yhat);
            \draw[arr] (Yhat) to (S);
        \end{tikzpicture}
        \caption{}
        \label{fig:pmar:selection:2}
    \end{subfigure}
    \begin{subfigure}[t]{0.32\linewidth}
        \centering
        \begin{tikzpicture}
            \node[var] (Y) at (1.3, 0) {$Y$};
            \node[var] (X) at (0, 0) {$X$};
            \node[var] (Z) at (.65, -1.2) {$Z$};
            \node[var] (S) at (.65, -2.4) {$S$};
            \draw[arr] (Y) to (X);
            \draw[biarr, bend left] (X) to (Y);
            \draw[biarr, bend left] (Y) to (Z);
            \draw[biarr, bend left] (Z) to (X);
            \draw[arr] (X) to (Z);
            \draw[arr] (Y) to (Z);
            \draw[arr] (Z) to (S);
        \end{tikzpicture}
        \caption{}
        \label{fig:pmar:selection:3}
    \end{subfigure}
    \caption{Causal graphs relating to examples of PMAR.}
\end{figure}

\paragraph{Examples of PMAR}
A causal graph related to the bank loan problem from the introduction is given in Figure \ref{fig:pmar:selection:1}, from which we deduce that indeed $Y\Indep S \given X, Z$. The bank loan problem is an example of a more general PMAR setting where individuals or items are \emph{selectively labelled} \citep{guerdan2023ground}, based on predictions from an existing model $\hat{Y} = \hat{\EE}[Y|X]$ and additional data $Z$ that is not part of the features $X$.

The PMAR setting is also encountered when we want to replace the feature set $Z$ of a current selection model $\hat{\EE}[Y|Z]$ with a new feature set $X$. One could start measuring $X$ when the old algorithm is still in use, to generate data on which the new prediction algorithm $\hat{\EE}[Y|X]$ can be learned; this data generating process is depicted in Figure \ref{fig:pmar:selection:2}.

Another example that fits the PMAR signature would be the selection of patients for costly CT scans, with $X$ digitally available measurements of the patient, $Z$ all additional information (besides $X$) that the doctor uses to decide whether a CT scan should be made of the patient (e.g. pale skin, slurred speech, pain indicated by patient), $Y$ whether the patient has a certain disease (as measured by the CT scan), and $S$ the doctors decision of making a CT scan of the patient. If we want to estimate $\EE[Y|X]$, e.g.\ to assist the doctor in the future in their decision making process, then we can only unbiasedly estimate this if we have measured $Z$, which requires the doctor to manually register all information that they use for making the decision. This can be a costly process, hence we might only want to have to measure this at train time, and not at test time.


The previous examples are of missing response, where the distributions $\PP(X, Y, Z  | S=1)$ and $\PP(X, Z)$ are readily available. As an example of selection bias with the PMAR conditional independence, we consider a dataset of patients from the University Hospital of Caracas, Venezuela, who have been tested for the presence of cervical cancer through a biopsy or colposcopy ($Y$).\footnote{Available at \url{https://archive.ics.uci.edu/ml/datasets/Cervical+cancer+(Risk+Factors)}.} In this dataset we have demographic and medical information of these patients ($X$), from which we might want to estimate $\EE[Y|X]$ for automated screening of the population. However, it might be that these patients are self-selected based on any symptoms, denoted by $Z$. A possible causal graph of the data generating process is depicted in Figure \ref{fig:pmar:selection:3}. Naive regression on this dataset would yield the biased model $\hat{\EE}[Y|X, S=1]$. As mentioned earlier, to correct for this bias we require external data $\Dcal$ from distribution $\PP(X,Z)$, i.e.\ an unbiased sample of covariates $X$ and symptoms $Z$ from the population. This could for example be sampled through a questionnaire. One could compare this with the approach of sampling $(X, Y)$ from the population to directly estimate $\EE[Y|X]$, which requires a costly and intrusive intervention (biopsy or colposcopy) on randomly sampled subjects, hence it would be preferred to only sample $(X, Z)$ from the population.

In all these cases, the variable $Z$ induces bias that should be corrected for when estimating $\EE[Y|X]$. In the following section, we demonstrate how this can be achieved.

\section{Imputation, Weighting, and Double Robustness under PMAR}\label{sec:methods}
\begin{figure}[!htb]
    \centering
    \fbox{\includegraphics[width=.95\linewidth]{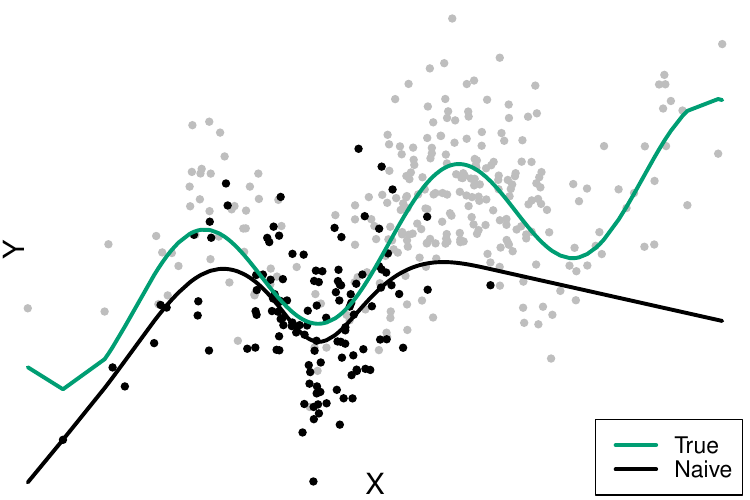}}
    \caption{Dataset from Example \ref{ex:1}, with observed (black) and unobserved data (grey); the true regression line $\EE[Y|X]$ (green) and the naively estimated regression line $\hat{\EE}[Y|X, S=1]$ (black).}
    \label{fig:ex1:true_and_naive}
\end{figure}

We propose three estimation procedures for training a regression model $\EE[Y|X]$ in the PMAR setting: a repeated regression method, a weighted regression method, and a doubly robust regression method. Throughout this section, we consider the following example.
\begin{example}\label{ex:1}
    Consider the data generating process
    \begin{align}\label{eqn:ex1}
        \begin{split}
            X &= \varepsilon_X \\
            Z &= 3\sin(X) + \varepsilon_Z \\
            Y &= \frac{1}{2}X + Z + \varepsilon_Y \\
            S &\sim \textrm{Bernoulli}(p_S(X, Z))
        \end{split}
    \end{align}
    with independent Gaussian random variables $\varepsilon_X, \varepsilon_Z$ and $\varepsilon_Y$. The selection probability is defined as $p_S(x, z) := \sigma(x)\sigma(z)$ with sigmoid $\sigma(x) = 1/(1+e^x)$.
\end{example}
\begin{wrapfigure}[7]{r}{0.28\linewidth}
    \centering
    \vspace{-15pt}
    \begin{tikzpicture}
        \node[var] (Y) at (1.3, 0) {$Y$};
        \node[var] (Z) at (1.3, -1.3) {$Z$};
        \node[var] (X) at (0, 0) {$X$};
        \node[var] (S) at (0, -1.3) {$S$};
        \draw[arr] (X) to (Z);
        \draw[arr] (X) to (Y);
        \draw[arr] (Z) to (Y);
        \draw[arr] (X) to (S);
        \draw[arr] (Z) to (S);
    \end{tikzpicture}
    \caption{}
    \label{fig:ex1:graph}
\end{wrapfigure}
Figure \ref{fig:ex1:graph} displays the graphical model related to this data generating process. Note that we can indeed read off that $Y \Indep S \given X, Z$, so if $Z$ is only available during training, the missingness mechanism is PMAR. Figure \ref{fig:ex1:true_and_naive} depicts $n=400$ draws from the generating model (\ref{eqn:ex1}). In this example, we have $\#\{S=1\}=123$. The black dots indicate $(X, Y)$ pairs where $Y$ is observed ($S=1$), and the grey dots indicate $(X, Y)$ pairs where $Y$ is unobserved ($S=0$). The green line shows the true regression line $\EE[Y|X=x] = \tfrac{1}{2}x + 3\sin(x)$, and the black line shows the regression line of the naive, biased estimate $\hat{\EE}[Y|X, S=1]$. In this example and throughout this paper we use thin plate splines regression \citep{duchon1977splines, wood2003thin} as implemented in the \href{https://cran.r-project.org/package=mgcv}{\code{mgcv}} package \citep{wood2015package}.\\
Note that, additional to what is shown in Figure \ref{fig:ex1:true_and_naive}, we have all $(X, Z)$ pairs available. Our goal is to fit a regression model that is close to the green line, $\EE[Y|X]$. \exqed

\paragraph{Extrapolation and positivity}
Many missingness methods assume \emph{positivity} of the selection probability, i.e.\ that $\PP(S=1 | X, Z) > 0$ holds $\PP(X, Z)$-almost surely, as is for example required for the identification equation (\ref{eqn:s-recoverability}). This assumption ensures that asymptotically, we do not run into the problem of extrapolation. For finite samples, we might run into extrapolation issues already when $\PP(S=1 | X, Z)$ is close to zero. In the example of Figure \ref{fig:ex1:true_and_naive} we have no positivity for large values of $X$; we will see that extrapolation is still possible by using privileged data $Z$ that is highly predictive of $Y$.

\subsection{Repeated regression}\label{subsec:rr}
By the law of total expectation and the PMAR conditional independence $Y\Indep S \given X, Z$ we can write the conditional expectation $\EE[Y|X]$ as follows
\begin{align}
    \begin{split}
        \EE[Y|X] &= \EE[\EE[Y|X, Z] | X] \\
        &= \EE[\EE[Y|X, Z, S=1] | X],
    \end{split}
\end{align}
similar to equation (\ref{eqn:s-recoverability}). For this equality to hold we require $\PP(S=1 | X, Z) > 0, \PP(X, Z)$-almost surely. We formulate an estimation procedure based on this expression by estimating each conditional expectation with a regression model. That is, we first regress $Y$ on $X$ and $Z$ using the data $\Scal$, which results in the regression model $\hat{\EE}[Y|X, Z, S=1]$. Using the unbiased data $(x_i, z_i)$ for $i\in\Dcal$ we construct \emph{pseudo-labels} $\tilde{Y}_i := \hat{\EE}[Y|X=x_i, Z_i=z_i, S=1]$. Now, we regress $\tilde{Y}$ on $X$ using the data $(X_i, \tilde{Y}_i)_{i\in\Dcal}$, which yields the final model:
\begin{equation}\label{eqn:rr}
    \hat{\mu}_{RR}(x) = \hat{\EE}[\tilde{Y}|X = x].
\end{equation}
Note that this method requires datasets $\Scal$ and $\Dcal$, and assumes that $Y\Indep S \given X, Z$. So, if this conditional independence is satisfied, it can be directly applied in both the missingness setting, and in the selection bias setting where we have `external data' $\PP(X, Z)$.

\setcounter{example}{0}
\begin{example}[\textit{continued}]\label{ex:1:imputed}
    Recall that the data for $Y$ is generated by $Y = \frac{1}{2}X + Z + \varepsilon_Y$, so for fitting a model $\hat{\EE}[Y|X, Z, S=1]$, although we can only fit this model on a small part of the data (123 out of 400 observations) it is a relatively easy (linear) model to fit. We compute the pseudo-labels $\tilde{Y}_i := \hat{\EE}[Y|X = x_i, Z = z_i, S=1]$, as depicted with orange crosses in Figure \ref{fig:ex1:imputed_gam}. Regressing these imputed values on $X$ yields the final model $\hat{\mu}_{RR}$ (orange line in Figure \ref{fig:ex1:imputed_gam}). \exqed
    \begin{figure}[!htb]
        \centering
        \fbox{\includegraphics[width=.95\linewidth]{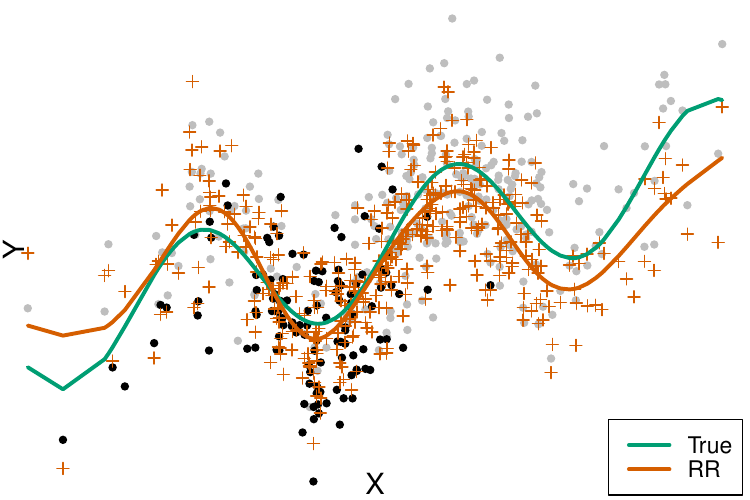}}
        \caption{Repeated regression. Orange crosses indicate imputed values, and the orange line is the regression line resulting from repeated regression.}
        \label{fig:ex1:imputed_gam}
    \end{figure}
\end{example}

Regression imputation is explicitly used for estimating the mean under missingness via $\EE[Y] = \EE[\EE[Y|X, S=1]]$ by \cite{bang2005doubly}. Similarly, the \emph{average causal effect} can be estimated by $\EE[Y|\Do(X)] = \EE[\EE[Y|X, Z]]$ for a valid \emph{adjustment set} $Z$; this method is known as \emph{standardization} \citep{hernan2021causal}. For these estimators, the inner expectation is estimated with a regression model, and the outer expectation is calculated by taking the mean of the pseudo-labels. The novelty in estimator (\ref{eqn:rr}) lies in the fact that the outer expectation is estimated by regression instead of the empirical mean, and the observation that this estimator is applicable to the PMAR setting.
Note that in general, automatically generated confidence intervals resulting from the outer regression $\hat{\EE}[\tilde{Y}|X]$ are not valid. To this end, \emph{multiple imputation} methods are often used. As this requires the modelling of the full distribution $\PP(Y|X, Z)$, we do not consider this. For an overview of multiple imputation methods in GEE regression, see \cite{ditlhong2018comparative}.

\subsection{Importance Weighting}\label{subsec:ipw}
\begin{figure}[!htb]
    \centering
    \fbox{\includegraphics[width=.95\linewidth]{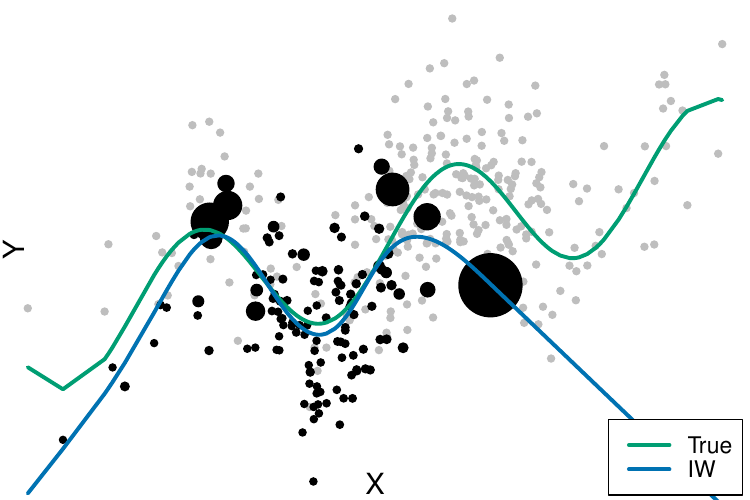}}
    \caption{IW regression with true weights.}
    \label{fig:ex1:ipw_known_weights}
\end{figure}

For estimating the parameter $\beta^*$ in the regression model $\EE[Y|X] = g(X; \beta^*)$ we often specify a loss function $\ell$ and perform empirical risk minimisation (\ref{eqn:erm}) as an approximation of the parameter that is optimal in terms of the true risk (\ref{eqn:risk}).
\begin{align}
    \hat{\beta} & = \argmin_\beta \frac{1}{|\Dcal|}\sum_{i\in\Dcal} \ell(g(X_i; \beta), Y_i) \label{eqn:erm} \\
    \beta^*     & = \argmin_\beta \EE[\ell(g(X; \beta), Y)] \label{eqn:risk}
\end{align}
Writing $f(x, y) := \ell(g(x; \beta), y)$, we can express the risk in terms of the distribution conditional on $S=1$ using importance weighting:
\begin{align}
    \begin{split}
        \EE[f(X, Y)]
        &= \EE[w(X, Z)f(X, Y) | S=1],
    \end{split}
\end{align}
with \emph{importance weights} $w(x, z) := \PP(S=1) / \PP(S=1 | X=x, Z=z)$, provided we have that $\PP(S=1 |X, Z) > 0$ holds $\PP(X, Z)$-almost surely. A derivation of these importance weights can be found in the supplementary material.
Since $\beta^* = \arg\min_\beta \EE[w(X, Z)\ell(g(X; \beta), Y) | S=1]$, when we have observations $(X_i, Z_i, Y_i) \sim \PP(X, Z, Y | S=1)$ for $i\in\Scal$, we can directly perform empirical risk minimization on this dataset using the weighted loss:
\begin{equation}
    \hat{\beta}_w := \argmin_\beta \frac{1}{|\Scal|}\sum_{i\in\Scal} w(X_i, Z_i)\ell(g(X_i; \beta), Y_i),
\end{equation}
and use $\hat{\beta}_w$ as an estimate of $\beta^*$. Many implementations of regression algorithms allow the user to specify such sample weights. When the weights are not known in the missingness setting, they can be estimated \citep{cortes2008sample}. Practically, when the selection probability $\PP(S=1 | x_i, z_i)$ is nearly zero for particular drawn values $x_i, z_i$, these probabilities are `clipped' (i.e.\ transformed to remain lower bounded by some pre-specified strictly positive value) to reduce variance.
In the experiments, we consider a linear map of the array of selection probabilities to $[1/20, 1]$, as this yields the best performance among different clipping strategies.
\setcounter{example}{0}
\begin{example}[\textit{continued}]\label{ex:1:ipw}
    The fitted IW regression model is depicted in Figure \ref{fig:ex1:ipw_known_weights}. Here, we used the true weights $w(x_i, z_i)$. The size of the observed points indicate the associated weight. Note that the IW regression only uses the black points (i.e.\ the $(x_i, y_i)$ pairs for which $S_i=1$).
    Comparing to repeated regression (Figure \ref{fig:ex1:imputed_gam}) we see that IW extrapolates poorly, possibly due to the effect of importance weighting on regularization \citep{xu2021understandinga}. On the other hand, IW interpolates better than the naive model. \exqed
\end{example}

The idea of weighting observations stems from the Horvitz-Thompson estimator for the population mean \citep{horvitz1952generalization}. The use of importance weights in GEE regression in the MAR setting is analysed by \citep{robins1994estimation,robins1995semiparametric}, and its MNAR counterpart is described by \cite{scharfstein1999adjusting}. When weighting is performed to correct for confounding bias in causal effect estimation, this is often called \emph{inverse propensity weighting} \citep{hernan2021causal}. For weighted SVM- and kernel regression under covariate shift and target shift, see \cite{zhang2013domain}.

\subsection{Doubly Robust Regression}\label{subsec:dr}
We now follow a relatively standard procedure to combine the two previous models into a doubly robust model \citep{bang2005doubly, kang2007demystifying, dudik2014doubly}. First, we consider the repeated regression method $\hat{\mu}_{RR}$ from Section \ref{subsec:rr}. We calculate the residuals of this method on the available $Y$ values: we set $R_i := Y_i - \hat{\mu}_{RR}(X_i)$ for all $i\in \Scal$. We model the conditional expectation of the residuals given $X$ using the IW regression method from Section \ref{subsec:ipw}, as
\begin{equation}
    \hat{r}_{IW}(x) := g(X; \hat{\beta}_w) \approx \EE[R | X]
\end{equation}
where $\hat{\beta}_w = \argmin_\beta |\Scal|^{-1}\sum_{i\in\Scal} w(X_i, Z_i)\ell(g(X_i; \beta), R_i)$ and where we use the same importance weights $w$ as in Section \ref{subsec:ipw}. Then we define the \emph{doubly robust} model
\begin{equation}
    \hat{\mu}_{DR}(x) := \hat{\mu}_{RR}(x) + \hat{r}_{IW}(x).
\end{equation}
Here, double robustness refers to the fact that for this model to be consistent, only one of $\hat{\mu}_{RR}$ and $\hat{r}_{IW}$ has to be consistent.

\setcounter{example}{0}
\begin{example}[\textit{continued}]\label{ex:1:dr}
    For sake of exposition, we apply the doubly robust method with a misspecified model $\hat{\mu}_{RR}(x)$. Recall that $\hat{\mu}_{RR}(x) = \hat{\EE}[\tilde{Y}|X= x]$, where $\tilde{Y}$ are the values that are imputed with the model $\hat{\EE}[Y|X, Z, S=1]$. To make the RR model deliberately misspecified, we use the same $\tilde{Y}$ as in Section \ref{subsec:rr} (more specifically, as in Figure \ref{fig:ex1:imputed_gam}), but for the outer regression $\hat{\EE}[\tilde{Y}|X]$ we employ polynomial regression with degree 5. The resulting models are depicted in Figure \ref{fig:ex1:dr}.
    We see that $\hat{\mu}_{DR}$ extrapolates poorly, but it interpolates better than the misspecified model $\hat{\mu}_{RR}$. \exqed
    \begin{figure}[!htb]
        \centering
        \fbox{\includegraphics[width=.95\linewidth]{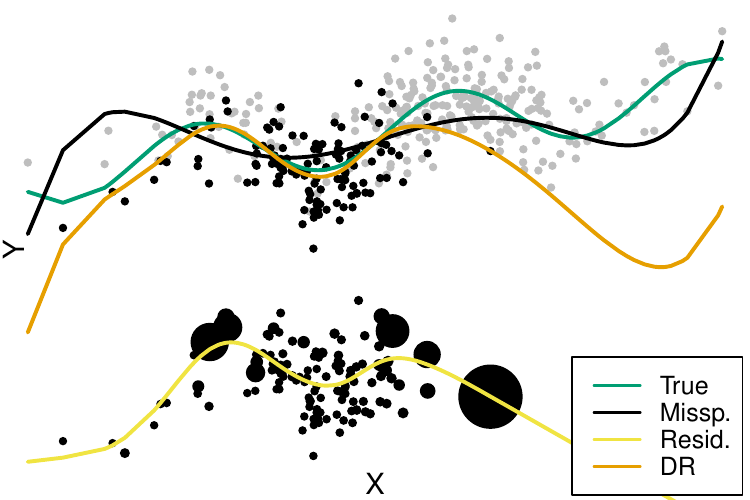}}
        \caption{Doubly robust regression. At the top the $(X,Y)$ pairs are plotted, with the true model (green) and misspecified RR model (black). The residuals of the misspecified model are plotted at the bottom. The IW regression model for the residuals (yellow) is then added to the misspecified model, which yields the doubly robust model (orange).}
        \label{fig:ex1:dr}
    \end{figure}
\end{example}
For overview papers on doubly robust estimation of the mean, see \cite{bang2005doubly} and \cite{kang2007demystifying}. For doubly robust regression in the MNAR setting, see \cite{rotnitzky1998semiparametric} and \cite{scharfstein1999adjusting}. Doubly robust estimation of causal effects is treated by \cite{coston2020counterfactual} and \cite{bhattacharya2022semiparametrica}, among others.

\section{Evaluation}\label{sec:evaluation}
Consider the setup where we split the available data $\Dcal = \Dcal'' \dot{\cup} \Dcal'$ into train- and test test respectively, and where we split $\Scal$ accordingly. When all values of $Y$ are known, performance is typically assessed using $\textrm{MSE} := |\Dcal'|^{-1}\sum_{i\in\Dcal'}(\hat{y}_i - y_i)^2$. Under missingness or selection bias, this quantity cannot be calculated as it requires the unobserved values $y_i$ for $i\notin \Scal'$. Quantities that \emph{can} be calculated on the test set are as follows:
\begin{align}
    \textrm{MSE-n}         & := \frac{1}{|\Scal'|}\sum_{i \in \Scal'} (\hat{y}_i - y_i)^2            \\
    \textrm{MSE-}w         & := \frac{1}{|\Scal'|}\sum_{i \in \Scal'} w(x_i, z_i)(\hat{y}_i - y_i)^2 \\
    \textrm{MSE-}\tilde{y} & := \frac{1}{|\Dcal'|}\sum_{i \in \Dcal'} (\hat{y}_i - \tilde{y}_i)^2,
\end{align}
which are respectively the naively calculated MSE, the re-weighted MSE using the true weights, and the MSE calculated as if the imputed values $\tilde{y}$ were true. We define $\textrm{MSE-}\hat{w}$ similarly, in terms of the estimated weights.
\setcounter{example}{0}
\begin{example}[\textit{continued}]\label{ex:1:eval}
    The IW and DR models are learned both with true and estimated weights (denoted with -t and -e respectively), and all weights are clipped. We calculate the different mean squared error metrics for the estimated models on 500 independently drawn test sets of $n=400$ samples. The results are provided in Table \ref{tab:ex1:eval}.
    \begin{table}[!htb]
        \centering
        \begin{tabular}{lccccc}\toprule
                  & $\textrm{MSE}$ & $\textrm{MSE-n}$ & $\textrm{MSE-}\tilde{y}$ & $\textrm{MSE-}w$ & $\textrm{MSE-}\hat{w}$ \\ \midrule
            Naive & 24.07          & \textbf{7.56}    & 19.54                    & 16.19            & 19.79                  \\
            RR    & \textbf{8.81}  & 8.13             & \textbf{4.44}            & \textbf{8.50}    & \textbf{9.34}          \\
            IW-t  & 40.42          & 8.40             & 35.58                    & 22.22            & 24.67                  \\
            IW-e  & 38.44          & 8.28             & 33.60                    & 21.49            & 23.08                  \\
            DR-t  & 25.91          & 8.13             & 21.27                    & 14.82            & 16.65                  \\
            DR-e  & 24.09          & 8.00             & 19.43                    & 13.96            & 15.52                  \\\midrule
            True  & 8.03           & 7.78             & 4.63                     & 7.87             & 8.87                   \\
            \bottomrule
        \end{tabular}
        \caption{Mean squared errors of different regression methods applied to Example \ref{ex:1:eval}, averaged over 500 test sets.}
        \label{tab:ex1:eval}
    \end{table}
    We emphasize that MSE is not calculable in practice, and note that MSE-n, MSE-$w$ and MSE-$\hat{w}$ yield qualitatively different results than the oracle MSE. MSE-$\tilde{y}$ provides the same ordering as MSE, but in practice the reliance on the imputation model can be undesirable. We observe that RR appropriately corrects for bias, as its MSE lies close to the MSE of the true function $\EE[Y|X]$. Additional to the oracle MSE, we calculate the MSEs specifically for data points where the regression model is interpolating (between the minimum and maximum values of $X$ for which $S=1$) and extrapolating (the complement of these values), as depicted in Table \ref{tab:example1_interp_extrap}. We observe that IW extrapolates poorly, especially compared to RR. \exqed
\end{example}
\begin{table}[!htb]
    \centering
    \begin{tabular}{lccccc}\toprule
             & $\textrm{MSE}$ & $\textrm{MSE-interp.}$ & $\textrm{MSE-extrap.}$ \\ \midrule
        RR   & \textbf{8.81}  & \textbf{8.43}          & \textbf{10.17}         \\
        IW-t & 40.42          & 12.00                  & 138.22                 \\
        IW-e & 38.44          & 11.77                  & 129.41                 \\
        DR-t & 25.91          & 10.47                  & 78.48                  \\
        DR-e & 24.09          & 10.23                  & 70.42                  \\
        \bottomrule
    \end{tabular}
    \caption{Inter- and extrapolation results for Example \ref{ex:1:eval}.}
    \label{tab:example1_interp_extrap}
\end{table}

When the true weights are not known, one can only evaluate the proposed bias correction methods on biased data by relying on an auxiliary model, i.e. either the imputation model $\hat{\EE}[Y|X, Z, S=1]$ or $\hat{\PP}(S=1 | X, Z)$. In general, if the practitioner is able to evaluate on an unbiased random sample of $(X, Y)$, this would be much more reliable than evaluating on a biased sample.

\section{Experiments}\label{sec:experiments}
The example that is used throughout this paper is quite extreme in the sense that the model $\EE[Y|X, Z]$ is relatively simple and $\EE[Y|X]$ is rather complex, which gives RR a clear advantage over IW. To investigate the proposed methods in a more general setting, we assess performance on extensive simulations and on the Boston Housing dataset.\footnote{Code for the experiments is publicly available at \url{https://github.com/philipboeken/debiased_regression}.}

\subsection{Simulations}\label{subsec:exp:sim}
To assess the performance of the proposed methods empirically, we first identify which graphs with variables $X, Y, Z$ and $S$ satisfy the pattern of d-separations $X\nPerp Y$ (such that regression $\EE[Y|X]$ makes sense), $Y\nPerp S|X$ and $Y\Perp S\given X, Z$ (PMAR). There are 550 \emph{Acyclic Directed Mixed Graphs} (ADMGs) that fit this pattern.
For each ADMG we simulate 50 datasets, and for each dataset we draw $2000$ samples according to an additive noise structural equation model, where all equations for $X, Y, Z$ are random functions, independently drawn from a Gaussian process, and the additive noise is the pushforward of a $\textrm{Unif}[0,1]$ distribution through a random function that is also drawn from a Gaussian process \citep{mooij2016distinguishing}. We let the selection probability be $p(x, z) = \sigma(x)\sigma(z)$, when $X$ and $Z$ are parents of $S$ in the ADMG. We make a 50/50 train-test split for training and evaluation. For a complete description of the simulation setup, we refer to the supplements.
On each dataset, we fit a `naive' model $\hat{\EE}[Y|X, S=1]$ and the models $\hat{\mu}_{RR}, \hat{\mu}_{IW}$ and $\hat{\mu}_{DR}$ from Section \ref{sec:methods}. We also fit a `true' model $\hat{\EE}[Y|X]$ for comparison.
\begin{table}[!htb]
    \centering
    \begin{tabular}{lcccc}\toprule
              & $\textrm{MSE}$               & $\textrm{MSE-}\tilde{y}$     & $\textrm{MSE-}w$             & $\textrm{MSE-}\hat{w}$       \\ \midrule
        Naive & 3.11 {\small (20.6)}         & 2.81 {\small (18.0)}         & 0.91 {\small (1.3)}          & 0.90 {\small (0.6)}          \\
        RR    & \textbf{2.01} {\small (2.6)} & \textbf{0.60} {\small (0.9)} & 1.25 {\small (1.5)}          & 1.20 {\small (1.3)}          \\
        IW-t  & 4.91 {\small (33.6)}         & 4.57 {\small (30.9)}         & 0.92 {\small (1.5)}          & 0.92 {\small (0.6)}          \\
        IW-e  & 4.18 {\small (23.2)}         & 3.83 {\small (20.7)}         & 0.92 {\small (1.2)}          & 0.91 {\small (0.6)}          \\
        DR-t  & 4.98 {\small (29.8)}         & 4.97 {\small (31.1)}         & \textbf{0.77} {\small (0.4)} & 0.81 {\small (0.4)}          \\
        DR-e  & 4.51 {\small (45.1)}         & 4.47 {\small (49.0)}         & 0.80 {\small (0.4)}          & \textbf{0.79} {\small (0.4)} \\ \midrule
        True  & 0.98 {\small (0.3)}          & 1.63 {\small (3.0)}          & 1.03 {\small (0.6)}          & 1.00 {\small (0.5)}          \\  \bottomrule
    \end{tabular}
    \caption{Results over 27.500 simulated datasets.}
    \label{tab:exp1_pos_indep}
\end{table}
In Table \ref{tab:exp1_pos_indep} we report the different mean squared error metrics of the methods, averaged over all 27.500 tests sets (standard deviations are shown in parentheses). We see that the average MSE of RR is half that of the best IW method. The median MSEs follow the relation $\textrm{RR} < \textrm{Naive} < \textrm{IW-t}$, confirmed with respective p-values $1.6\cdot 10^{-110}$ and $3.2\cdot 10^{-39}$ of Wilcoxon's signed rank test. In ADMGs where $X$ is a parent of $S$ there is a clear region of the support of $X$ with no positivity of $S$, hence a clear region where the regression model is extrapolating, similar to Example \ref{ex:1}. When singling out these ADMGs we calculate the average MSE, and additionally the MSEs on the parts where the regression model is inter- and extrapolating, as shown in the table below. From Table \ref{tab:exp1_interp_extrap} we see that IW and DR do not extrapolate well, while RR does.

\begin{table}[!htb]
    \centering
    \begin{tabular}{lccccc}\toprule
             & $\textrm{MSE}$               & $\textrm{MSE-interp.}$       & $\textrm{MSE-extrap.}$       \\ \midrule
        RR   & \textbf{2.13} {\small (2.7)} & 1.46 {\small (1.8)}          & \textbf{2.89} {\small (4.2)} \\
        IW-t & 7.29 {\small (34.0)}         & 1.29 {\small (0.7)}          & 13.75 {\small (53.8)}        \\
        IW-e & 5.82 {\small (16.4)}         & 1.27 {\small (0.6)} & 10.81 {\small (32.6)}        \\
        DR-t & 7.73 {\small (45.7)}         & 1.26 {\small (0.6)}          & 14.37 {\small (70.1)}        \\
        DR-e & 6.93 {\small (72.3)}         & \textbf{1.24} {\small (0.6)} & 12.34 {\small (94.9)}        \\
        \bottomrule
    \end{tabular}
    \caption{Inter- and extrapolation results for simulated data.}
    \label{tab:exp1_interp_extrap}
\end{table}

We hypothesize that IW and DR extrapolate badly because of the interplay between the smoothness assumption and regularization of thin plate spline regression, and large weights near the edge of the support of $\PP(X | S=1)$, causing the model to diverge away from the true $\EE[Y|X]$ when extrapolating (as is for example the case in Figure \ref{fig:ex1:ipw_known_weights}). We investigate this by running the experiments with RR, IW and DR implemented with regression trees \citep{breiman1984classification} instead of thin plate regression. As regression trees extrapolate flatly, this can allow for better extrapolation performance than thin pate regression. 
The results of these experiments are provided in the supplements. They are indeed numerically less extreme than for thin plate regression, but qualitatively the same.

\subsection{Boston Housing data}\label{subsec:exp:rw}
In absence of a real-world dataset where there is missingness \emph{and} the missing value is known,
we consider the Boston Housing dataset \citep{harrison1978hedonic} and simulate a selection mechanism ourselves.
We consider the problem of predicting the value of owner-occupied homes in US Dollars ($Y$) from the number of rooms per dwelling ($X$). We let selection depend on $X$, and on the `percentage of people of lower status of the population in the town where the house is situated' ($Z$), which strongly correlates with $Y$. Selection is simulated by setting $\PP(S=1|X, Z) := \sigma(f_1(X))\sigma(f_2(Z))$, where $f_1, f_2$ are random functions drawn from a Gaussian process. For each dataset, the 506 available observations are randomly split into evenly sized train and test sets. The average MSEs are provided in Table \ref{tab:exp2}, and show the same qualitative results as the simulations.
When calculating the MSE of the inter- and extrapolation regions, we observe that RR performs better than IW at both tasks (table shown in supplements).
\begin{table}[!htb]
    \centering
    \begin{tabular}{lccccc}\toprule
              & $\textrm{MSE}$               & $\textrm{MSE-}\tilde{y}$     & $\textrm{MSE-}w$             & $\textrm{MSE-}\hat{w}$       \\ \midrule
        Naive & 1.23 {\small (2.5)}          & 0.84 {\small (2.3)}          & 0.48 {\small (0.7)}          & 1.79 {\small (5.9)}          \\
        RR    & \textbf{0.71} {\small (0.3)} & \textbf{0.26} {\small (0.3)} & 0.44 {\small (0.4)}          & 1.63 {\small (4.8)}          \\
        IW-t  & 2.18 {\small (4.9)}          & 1.77 {\small (4.9)}          & 0.56 {\small (1.5)}          & 2.40 {\small (9.2)}          \\
        IW-e  & 1.75 {\small (4.4)}          & 1.40 {\small (4.7)}          & 0.50 {\small (0.9)}          & 1.98 {\small (7.4)}          \\
        DR-t  & 1.92 {\small (3.7)}          & 1.65 {\small (3.4)}          & \textbf{0.23} {\small (0.2)} & 0.39 {\small (1.2)}          \\
        DR-e  & 2.43 {\small (5.4)}          & 2.20 {\small (5.5)}          & 0.25 {\small (0.2)}          & \textbf{0.21} {\small (0.3)} \\\midrule
        True  & 0.45 {\small (0.1)}          & 0.54 {\small (0.4)}          & 0.36 {\small (0.3)}          & 0.85 {\small (2.2)}          \\
        \bottomrule
    \end{tabular}
    \caption{Results of $500$ biased instantiations of the Boston Housing dataset.}
    \label{tab:exp2}
\end{table}

\section{Discussion and conclusion}
In this work, we have motivated the use of \emph{privileged information} for estimating a regression model when selection or missingness is nonignorable, and introduced the \emph{Privilegedly Missing at Random} (PMAR) setting. We formulated the repeated regression method, the IW regression method, and the doubly robust combination of the two. We note that evaluation of regression methods on biased data is not straightforward and relies on auxiliary models. Experiments show that repeated regression can appropriately correct for bias, and with considerable advantage over IW regression. In particular, in repeated regression, extrapolation is facilitated by privileged data $Z$ that is predictive of $Y$. IW does not have such a property, and extrapolates worse than RR.

Further research can be done on the statistical properties of the proposed methods, e.g.\ on bounding the regret of DR in terms of the regret of RR and IW (possibly under misspecification), or on the interplay of importance weighting, regularization, and extrapolation. 

Additionally, further research on the sensitivity of these methods with respect to the conditional independence assumption would be valuable, e.g.\ to determine whether it's always better to correct for bias with some privileged information $Z$, even when the independence $Y\Indep S\given X, Z$ is not met.

Finally, it would be interesting to see whether the improved performance of repeated regression over weighted regression translates to their counterparts for estimating causal effects, i.e.\ whether standardization should be preferred over inverse propensity weighting.

\begin{acknowledgements}
    This work is supported by Booking.com. We thank Kees Jan de Vries for fruitful discussions and Christina Katsimerou, Nils Skotara, Kostas Tokis and anonymous reviewers for their feedback on the manuscript.
\end{acknowledgements}

\bibliographystyle{abbrvnat}
\renewcommand{\bibsection}{\subsubsection*{References}}
\bibliography{philip}

\end{document}


\onecolumn

\maketitle

\section{Selection bias and Missingness}
A schematic display of the available data under missingness and under selection bias with external data is provided in Figure \ref{fig:sbias_vs_missingness}. Note that under missingness we can estimate $\PP(S=1|X, Z)$ from the dataset $\Dcal$; under selection bias this is not possible.

\begin{figure}[!htb]
    \centering
    \begin{subfigure}{\linewidth}
        \centering
        \begin{tabular}{ccllclc}
            \cmidrule[\heavyrulewidth]{3-6}
                                              &       & $X$                                      & $Z$       & $S$                                                      & $Y$                            & \\
            \cmidrule[\heavyrulewidth]{3-6}
            \ldelim\{{7}[-3pt]{8pt}[$\Dcal$]  &
            \ldelim\{{3}[-3pt]{10pt}[$\Scal$] &
            $x_1$                             & $z_1$ & 1                                        & $y_1$     & \multirow{3}{*}{\hspace{-8pt}\small$\PP(X, Y, Z | S=1)$}                                    \\
                                              &       & \vdots                                   &           &                                                          &                                & \\
                                              &       & $x_m$                                    & $z_m$     & 1                                                        & $y_m$                            \\
            \cmidrule{3-6}
                                              &       & $x_{m+1}$                                & $z_{m+1}$ & 0                                                        & $y_{m+1}$                      & \\
                                              &       & \vdots                                   &           &                                                          &                                & \\
                                              &       & $x_n$                                    & $z_n$     & 0                                                        & $y_n$~~~~~\boxit{23pt}{36.4pt} & \\
            \cmidrule[\heavyrulewidth]{3-6}
                                              &       & \multicolumn{3}{c}{\small$\PP(X, Z, S)$} &           &
        \end{tabular}
        \caption{Missing response}
    \end{subfigure}\\
    \vspace{20pt}
    \begin{subfigure}{.8\linewidth}
        \centering
        \begin{minipage}[t][][b]{.55\linewidth}
            \centering
            \begin{tabular}{cllclc}
                \cmidrule[\heavyrulewidth]{2-5}
                                                  & $X$       & $Z$       & $S$   & $Y$                                                      & \\
                \cmidrule[\heavyrulewidth]{2-5}

                \ldelim\{{3}[-3pt]{10pt}[$\Scal$] &
                $x_1$                             & $z_1$     & 1         & $y_1$ & \multirow{3}{*}{\hspace{-8pt}\small$\PP(X, Y, Z | S=1)$}   \\
                                                  & \vdots    &           &       &                                                          & \\
                                                  & $x_m$     & $z_m$     & 1     & $y_m$                                                      \\
                \cmidrule{2-5}
                                                  & $x_{m+1}$ & $z_{m+1}$ & 0     & $y_{m+1}$                                                & \\
                                                  & \vdots    &           &       &                                                          & \\
                                                  & $x_n$     & $z_n$     & 0     & $y_n$~~~~~\boxit{113pt}{36.4pt}                          & \\
                \cmidrule[\heavyrulewidth]{2-5}
            \end{tabular}
        \end{minipage}%
        \begin{minipage}[t][][b]{.35\linewidth}
            \centering
            \begin{tabular}{cll}
                \cmidrule[\heavyrulewidth]{2-3}
                                                    & $X$                                   & $Z$    \\
                \cmidrule[\heavyrulewidth]{2-3}
                \ldelim\{{6.3}[-3pt]{10pt}[$\Dcal$] & $x_1$                                 & $z_1$  \\
                                                    & \vdots                                & \vdots \\
                                                    & \vdots                                & \vdots \\
                                                    & \vdots                                & \vdots \\
                                                    & $x_n$                                 & $z_n$  \\
                \cmidrule[\heavyrulewidth]{2-3}
                                                    & \multicolumn{2}{c}{\small$\PP(X, Z)$}
            \end{tabular}
        \end{minipage}
        \caption{Selection bias with external data}
    \end{subfigure}%
    \caption{Available data under missingness and under selection bias with external data. Grayed-out areas indicate unobserved data.}
    \label{fig:sbias_vs_missingness}
\end{figure}

In the PMAR setting we have $Y\Indep S \given X, Z$ and we are only given values of $X$ at test time; our target is to estimate the function $\EE[Y|X]$.

\newpage
\section{Importance weighting}
For estimating the parameter $\beta^*$ in the regression model $\EE[Y|X] = g(X; \beta^*)$ we often specify a loss function $\ell$ and perform empirical risk minimisation (\ref{eqn:erm}) as an approximation of the optimal parameter in terms of the true risk (\ref{eqn:risk}).
\begin{align}
    \hat{\beta} & = \argmin_\beta \frac{1}{n}\sum_{i=1}^n \ell(g(X_i; \beta), Y_i) \label{eqn:erm} \\
    \beta^*     & = \argmin_\beta \EE[\ell(g(X; \beta), Y)] \label{eqn:risk}
\end{align}
Writing $f(x, y) := \ell(g(x; \beta), y)$, we can express the risk in terms of the distribution conditional on $S=1$ using importance weighting:
\begin{align}
    \begin{split}
        \EE[f(X, Y)] &= \int f(x, y)p(x, y, z)\diff (x, y, z) \\
        &= \int f(x, y)p(x, y, z)\frac{p(x, y, z | S=1)}{p(x, y, z | S=1)}\diff (x, y, z) \\
        &= \int f(x, y)\frac{p(x, y, z)\PP(S=1)}{p(x, y, z, S=1)}p(x, y, z | S=1)\diff (x, y, z) \\
        &= \int f(x, y)\frac{\PP(S=1)}{\PP(S=1 | x, y, z)}p(x, y, z | S=1)\diff (x, y, z) \\
        &= \int f(x, y)\frac{\PP(S=1)}{\PP(S=1 | x, z)}p(x, y, z | S=1)\diff (x, y, z) \\
        &= \int f(x, y)w(x, z)p(x, y, z | S=1)\diff (x, y, z) \\
        &= \EE[w(X, Z)f(X, Y) | S=1],
    \end{split}
\end{align}
where in the fifth equation we use the conditional independence $Y\Indep S \given X, Z$, and where we define the \emph{importance weights}
\begin{equation}
    w(x, z) := \frac{\PP(S=1)}{\PP(S=1 | x, z)}.
\end{equation}
Since $\beta^* = \arg\min_\beta \EE[w(X, Z)\ell(g(X; \beta), Y) | S=1]$, when we have observations $(X_1, Z_1, Y_1), ..., (X_n, Z_n, Y_n) \sim \PP(X, Z, Y | S=1)$, we can directly perform empirical risk minimization on this dataset using the weighted loss:
\begin{equation}
    \hat{\beta} = \argmin_\beta \frac{1}{n}\sum_{i=1}^n w(x_i, z_i)\ell(g(x_i; \beta), y_i).
\end{equation}

\section{Simulations}
\begin{figure}[!htb]
    \centering
    \includegraphics[width=\textwidth]{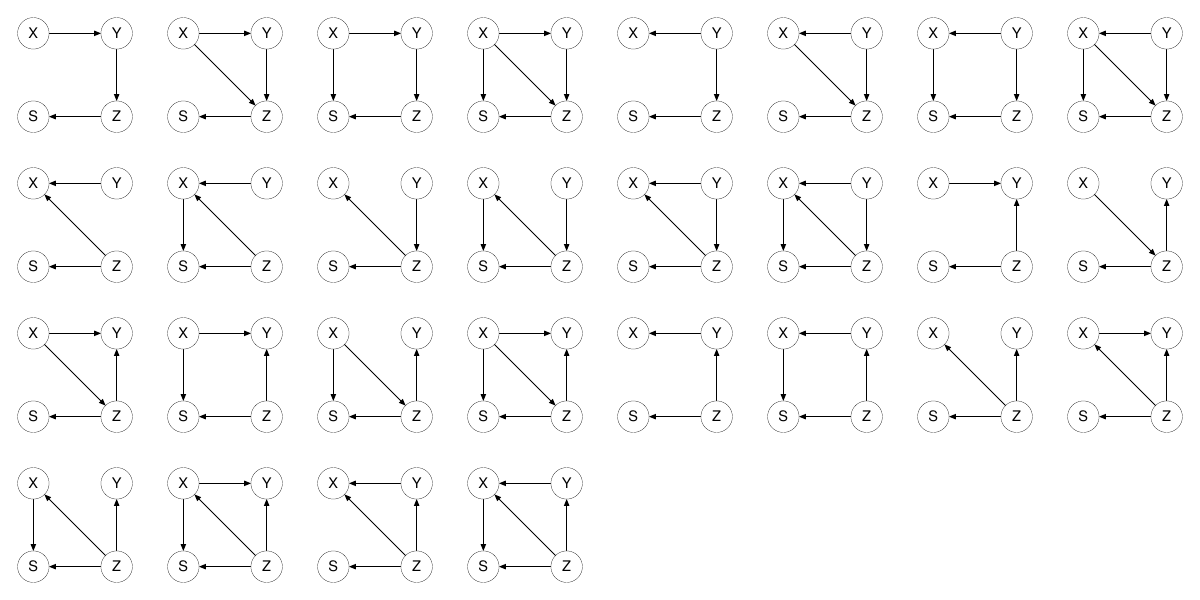}
    \caption{PMAR DAGs where $S$ is a sink node.}
    \label{fig:dag_s_sink}
\end{figure}

\begin{figure}[!htb]
    \centering
    \includegraphics[width=\textwidth]{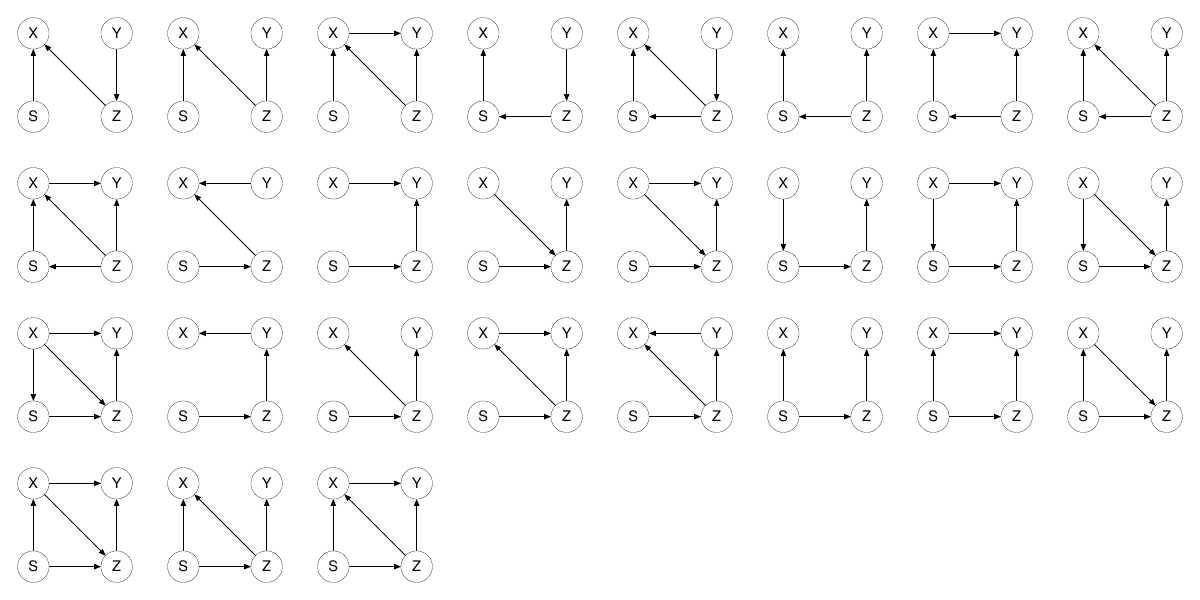}
    \caption{PMAR DAGs where $S$ is not a sink node.}
    \label{fig:dag_s_not_sink}
\end{figure}

We performed a brute-force search of all \emph{Acyclic Directed Mixed Graphs} (ADMGs) that satisfy $X\nPerp Y$ and the PMAR pattern of d-separations $Y\nPerp S$, $Y\nPerp S\given X$, and $Y\Perp S \given X, Z$, using the \code{pcalg} package \citep{kalisch2012causal}, resulting in 550 ADMGs. For reference, all 55 DAGs are depicted in figures \ref{fig:dag_s_sink} and \ref{fig:dag_s_not_sink}, where the graphs are categorised by whether $S$ has any children or not. The remaining 495 ADMGs are not depicted here.

For each of the 550 ADMGs, we simulate 50 datasets according to the following procedure. Throughout, let $V, S, X$ etc.\ denote vertices in the graph, and let $x_V, x_S, x_X$ denote vectors of dimension $n=2000$ for the simulated values.
\begin{itemize}
    \item First, we replace any bidirected edge pointing to variables $\bm{V}$ by a variable $U_{\bm{V}}$, and let $\bm{V}$ be the children of $U_{\bm{V}}$. This turns the ADMG into a DAG.
    \item Then, we calculate a topological order of the DAG.
    \item For every variable $V$ in the topological order, we simulate 2000 observations as follows:
    \begin{itemize}
        \item If $V$ has no parents, then
        \begin{itemize}
            \item if $V=S$, draw $x_S \sim \textrm{Bernoulli}(1/3)$;
            \item otherwise, $V\neq S$ and we draw $x_V \sim \Rcal\Dcal$, as defined below.
        \end{itemize}
        \item Otherwise, denote the parents of $V$ with $\bm{Pa}$ and their value $\bm{x}_{\bm{Pa}}$, and then
        \begin{itemize}
            \item if $V = S$, draw $(x_S)_i \sim \textrm{Bernoulli}(p((x_{\bm{Pa}})_i))$ where
            \begin{equation}
                p((x_{\bm{Pa}})_i) := \prod_{v\in \bm{Pa}}\sigma((x_v)_i)
            \end{equation}
            and $\sigma(x) := (1+e^{20x})^{-1}$;
            \item otherwise, $V\neq S$ and
            \begin{itemize}
                \item draw a random function $f_V$ from a Gaussian process on $\RR^{|\bm{Pa}\setminus \{S\}|}$ as $f_V\sim\Gcal\Pcal(0, K_M)$;
                \item draw noise $\varepsilon_V \sim \Rcal\Dcal$ and set
                \begin{equation}
                    x_V := f_V(\bm{x}_{\bm{Pa}\setminus \{S\}}) + \frac{1}{2}\varepsilon_V;
                \end{equation}
                \item if $S\in \bm{Pa}$, calculate the empirical standard deviation $c := \textrm{sd}(x_V)$ and set $(x_V)_i := (x_V)_i - 3c$ for all $i$ where $(x_S)_i = 1$;
                \item then, standardize $x_V$.
            \end{itemize}
        \end{itemize}
    \end{itemize}
\end{itemize}

\paragraph{Drawing from a Gaussian process} The kernels used for calculating the covariance matrix for drawing from a Gaussian process are the Mat\'ern kernel and the squared exponential kernel, being
\begin{align}
    K_M(x, y)     & := (1 + \sqrt{5} \|x-y\| + \frac{5}{3} \|x-y\|^2) e^{-\sqrt{5}} \label{eqn:matern} \\
    K_{SE} (x, y) & := \frac{1}{4} e^{\frac{2}{9}\|x-y\|^2}
\end{align}
respectively, where $\|\cdot\|$ denotes the Euclidean norm. Then, given input $x_i\in\RR^d$ for $i=1, ..., n$ and kernel $K:\RR^d\times \RR^{d} \to \RR$, denote $\bm{x} := (x_i)_{i=1}^n$ and draw $f(\bm{x}) \sim \Ncal(0, (K(x_i, x_j))_{i, j})$, where $i$ and $j$ run over $\{1, ..., n\}$ in the kernel matrix $(K(x_i, x_j))_{i, j}$. This draw is denoted with $f\sim \Gcal\Pcal(0, K)$.

\paragraph{Drawing from a random distribution} Drawing noise from a random distribution, denoted with $\varepsilon \sim \Rcal\Dcal$, is done as follows:
\begin{itemize}
    \item First, draw 2000 i.i.d.\ samples $U \sim \textrm{Unif}[0,1]$.
    \item Then, draw a random function $f_\varepsilon \sim \Gcal\Pcal(0, K_{SE})$
    \item Set $\varepsilon := f_\varepsilon(U)$, and standardize.
\end{itemize}

\clearpage

\subsection{Simulation experiments with regression trees}
In the main paper, we hypothesize that the bad extrapolation performance of weighted regression is caused by a detrimental effect of the weights on the regularization of the regression method. In the main paper we use thin plate spline regression, which doesn't necessarily extrapolate flatly and can diverge away from the true $\EE[Y|X]$, yielding large MSE values for IW and DR. Tables \ref{tab:exp1:tree:pos_indep} and \ref{tab:exp1:tree:interp_extrap} show the simulation results when using regression trees as implemented in the \href{https://cran.r-project.org/web/packages/rpart}{\code{rpart}} package, instead of thin plate regression. The results are numerically less extreme, but qualitatively the same as for thin plate regression, as depicted in the main paper.

\begin{table}[!htb]
\parbox[][145pt][t]{.5\linewidth}{%
    \begin{tabular}{lcccc}\toprule
            & $\textrm{MSE}$               & $\textrm{MSE-}\tilde{y}$     & $\textrm{MSE-}w$             & $\textrm{MSE-}\hat{w}$       \\ \midrule
        Naive & 1.53 {\small (0.6)}          & 0.62 {\small (0.8)}          & 0.96 {\small (0.6)}          & 0.96 {\small (0.5)}          \\
        RR    & \textbf{1.38} {\small (0.6)} & \textbf{0.27} {\small (0.3)} & 0.97 {\small (0.6)}          & 0.94 {\small (0.5)}          \\
        IW-t  & 1.54 {\small (0.6)}          & 0.63 {\small (0.8)}          & 0.97 {\small (0.6)}          & 0.98 {\small (0.5)}          \\
        IW-e  & 1.52 {\small (0.6)}          & 0.61 {\small (0.7)}          & 0.97 {\small (0.6)}          & 0.97 {\small (0.5)}          \\
        DR-t  & 1.48 {\small (0.7)}          & 0.58 {\small (0.7)}          & \textbf{0.66} {\small (0.3)} & 0.73 {\small (0.4)}          \\
        DR-e  & 1.46 {\small (0.6)}          & 0.55 {\small (0.7)}          & 0.71 {\small (0.4)}          & \textbf{0.69} {\small (0.3)} \\ \midrule
        True  & 1.00 {\small (0.2)}          & 0.69 {\small (0.7)}          & 1.05 {\small (0.6)}          & 1.02 {\small (0.5)}          \\   \bottomrule
    \end{tabular}
    \caption{Results over 27.500 simulated datasets.}
    \label{tab:exp1:tree:pos_indep}}
\parbox[][145pt][t]{.45\linewidth}{%
    \begin{tabular}{lccccc}\toprule
            & $\textrm{MSE}$               & $\textrm{MSE-interp.}$       & $\textrm{MSE-extrap.}$       \\ \midrule
        Naive & 1.58 {\small (0.6)}          & 1.37 {\small (0.6)}          & 1.82 {\small (1.0)}          \\
        RR    & \textbf{1.44} {\small (0.6)} & \textbf{1.22} {\small (0.6)} & \textbf{1.68} {\small (1.0)} \\
        IW-t  & 1.60 {\small (0.6)}          & 1.38 {\small (0.6)}          & 1.85 {\small (1.1)}          \\
        IW-e  & 1.58 {\small (0.6)}          & 1.37 {\small (0.6)}          & 1.82 {\small (1.0)}          \\
        DR-t  & 1.55 {\small (0.7)}          & 1.30 {\small (0.6)}          & 1.84 {\small (1.1)}          \\
        DR-e  & 1.52 {\small (0.6)}          & 1.28 {\small (0.6)}          & 1.80 {\small (1.0)}          \\ \midrule
        True  &  1.01 {\small (0.2)} & 1.03 {\small (0.3)} & 0.98 {\small (0.3)} \\ 
        \bottomrule
    \end{tabular}
    \caption{Interpolation and extrapolation results of regression trees for simulated data, on graphs with $X\to S$.}
    \label{tab:exp1:tree:interp_extrap}}
\end{table}

\section{Boston Housing Data}
\begin{figure}
    \centering
    \includegraphics[width=.65\linewidth]{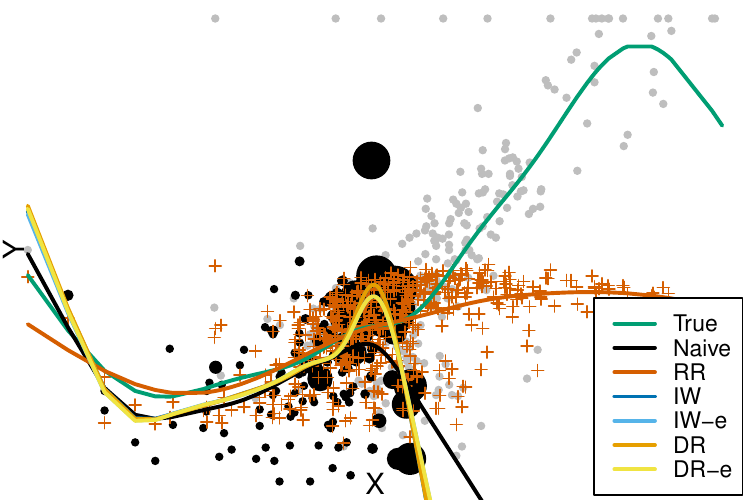}
    \caption{An instantiation of the biased Boston Housing dataset.}
    \label{fig:exp2}
\end{figure}

Considering the Boston Housing Dataset \citep{harrison1978hedonic}, let the variables $X, Y$ and $Z$ be `the number of rooms per dwelling', `the value of owner-occupied homes in US Dollars' and `percentage of people of lower status of the population' respectively. We sample the selection probability
\begin{equation}
    p(X, Z) := \sigma(f_1(X))\sigma(f_2(Z))
\end{equation}
where $f_1, f_2 \sim \Gcal\Pcal(0, K_{SE})$ independent. Then we draw $U_i\sim \textrm{Unif}[0,1]$ for $i=1, ..., 506$, and set
\begin{equation}
    S_i := \I\{p(X, Z)_i < U_i\}
\end{equation}
for all $i$. The dataset is resampled when $\#\{S = 1\} < 120$. One realisation of such a dataset with an overview of all regression methods is provided in Figure \ref{fig:exp2}.

The MSE and the interpolated and extrapolated variants, as calculated on the Boston Housing dataset, are shown in Table \ref{tab:boston_interp_extrap}. We observe that RR performs better than IW and DR on all three metrics.

\bibliographystyle{abbrvnat}
\renewcommand{\bibsection}{\subsubsection*{References}}
\bibliography{philip}

\begin{table}
    \centering
    \begin{tabular}{lccccc}\toprule
              & $\textrm{MSE}$               & $\textrm{MSE-interp.}$       & $\textrm{MSE-extrap.}$       \\ \midrule

        Naive & 1.23 {\small (2.5)}          & 0.57 {\small (0.5)}          & 3.88 {\small (9.1)}          \\
        RR    & \textbf{0.71} {\small (0.3)} & \textbf{0.47} {\small (0.2)} & \textbf{2.05} {\small (4.3)} \\
        IW-t  & 2.18 {\small (4.9)}          & 0.62 {\small (0.8)}          & 7.68 {\small (17.7)}         \\
        IW-e  & 1.75 {\small (4.4)}          & 0.58 {\small (0.6)}          & 5.96 {\small (16.1)}         \\
        DR-t  & 1.92 {\small (3.7)}          & 0.48 {\small (0.3)}          & 7.37 {\small (16.4)}         \\
        DR-e  & 2.43 {\small (5.4)}          & 0.52 {\small (0.4)}          & 9.42 {\small (22.1)}         \\
        \bottomrule
    \end{tabular}
    \caption{}
    \label{tab:boston_interp_extrap}
\end{table}